\documentclass[letterpaper, preprint, paper,11pt]{AAS}	% for preprint proceedings

\usepackage{bm}
\usepackage{amsmath}
\usepackage[colorlinks=true, pdfstartview=FitV, linkcolor=black, citecolor= black, urlcolor= black]{hyperref}
\usepackage{overcite}
\usepackage{footnpag}			      	% make footnote symbols restart on each page

\usepackage{amssymb}
\usepackage{multirow}
\usepackage{multicol}
\usepackage{subcaption}
\usepackage[dvipsnames]{xcolor}
\usepackage{ulem}
\PaperNumber{19-746}

\begin{document}

\title{A Study on Effective Initial Guess Finding Method Based on B\'ezier Curves: \ Orbit Determination Applications}

\author{Daegyun Choi\thanks{PhD Student, Department of Aerospace Engineering, Mississippi State University, Mississippi State, MS 39762, USA.}, 
Sungwook Yang\thanks{PhD Student, Department of Aerospace Engineering, Mississippi State University, Mississippi State, MS 39762, USA.},
Henzeh Leeghim\thanks{Associate Professor, Department of Aerospace Engineering, Chosun University, Gwangju 61452, Republic of Korea.},
\ and Donghoon Kim\thanks{Assistant Professor, Department of Aerospace Engineering \& Engineering Mechanics, University of Cincinnati, Cincinnati, OH 45221, USA.}
}

\maketitle{} 		

\begin{abstract}
In celestial mechanics, proper orbits related to missions are obtained by solving two-point boundary value problems. Since a selection method of initial value affects the convergence of the solution, developing an effective method to find an initial guess is required. In this work, B\'ezier curves, which can describe complicated curves and surfaces, are utilized to find the initial guess. First, the given problems are transformed into B\'ezier curves forms, and B\'ezier curves' control points, which can handle the shape of curves, are selected by solving the system of nonlinear equations. Finally, the initial guess is obtained by substituting the calculated control points to B\'ezier curves. To validate the performance of the proposed method, numerical simulations are conducted with respect to three kinds of orbits, which are from circular to highly elliptical orbit (HEO). The proposed method is compared to the general shooting method. The comparison results show that the initial guess calculated by B\'ezier curves makes finding the solution more efficient in terms of computational time and iterations. Also, it shows that the proposed method finds the solution for the HEO while the general shooting method fails to find the solution.
\end{abstract}

\section{Introduction}

In celestial mechanics, Lambert's problem is concerned with an orbit determination between an initial and final position with a time of flight (TOF). The orbit determination of spacecraft is usually performed to complete missions such as rendezvous, targeting, orbit transfer, etc. In general, determining the orbit is the one of the classical two-point boundary problems (TPBVPs) in celestial mechanics. To solve TPBVPs, the shooting method is generally used. The given problem is transformed into an initial value problem (IVP), and the solution of the TPBVP is obtained by solving the transformed IVP with a selection of initial values, which are not given. If the initial values are selected inappropriately, the solution of the TPBVP cannot converge to its boundary conditions. Therefore, developing an effective method to find its initial guess is required since the selection method of the initial values affects the convergence of the solution. For this reason, a new approach is applied to find an initial guess effectively. 

B\'ezier curves are parametric curves with great flexibility, which can describe complicated curves and surfaces\cite{Bezier1966,Bezier1967,Bezier1968} The functions of B\'ezier curves are linear combinations of polynomials with respect to control points, which play a role to handle the shape of curves. Because of the outstanding capability to describe the trajectory, B\'ezier curves have been investigated to find alternative approaches solving various problems in orbital mechanics. In 2007, Mortari and Clocchiatti\cite{Mortari2007} applied non-rational B\'ezier curves through approximating Kepler's equation for the entire mean anomaly range, and a few years later, B\'ezier curves were used as a root-finding tool for solving Kepler's equation\cite{Mortari2014}. Also, a formulation of Lambert's problem in terms of the argument of periapsis using B\'ezier curves was designed by Pan and Ma\cite{Pan2018}. For conic section descriptions, Kim and Mortari\cite{Kim2014_2} presented closed-form solutions for midpoint weights using rational B\'ezier curves, and methodologies to find analytical expressions for TPBVPs were proposed using non-rational B\'ezier curves\cite{Kim2014}. In this work, B\'ezier curves are utilized to finding an initial guess to solve TPBVPs, especially the orbit determination problem. Using the definition of B\'ezier curves, given TPBVPs are transformed into B\'ezier function forms, which are linear combinations of polynomials with control points. After the control points are selected properly, the initial guess is obtained by substituting the control points with variables of B\'ezier curves. Numerical simulations are performed to validate the performance of the method proposed with respect to a nonlinear TPBVP and Lambert's problem. The results using the initial guess obtained by B\'ezier curves are compared to the results using an arbitrary initial guess.

\section{B\'ezier Curves}
B\'ezier curves are based on Bernstein polynomials, which can describe complicated curves and surfaces\cite{Bezier1966,Bezier1967,Bezier1968}. B\'ezier curves are defined by a set of control points ${\bf P}_0$ though ${\bf P}_n$ ($n+1$ number of points), where $n$ is B\'ezier curves' order. The first and last control points on the curve are called endpoints, and they indicate where the curve starts and ends. The intermediate control points generally do not lie on the curve\cite{Farin2002}, but they control the shape of the curve. There are two types of B\'ezier curves, such as non-rational and rational B\'ezier curves. Rational B\'ezier curves require more variables to include adjustable weights than non-rational B\'ezier curves. In this work, non-rational B\'ezier curves are considered in order to reduce the number of unknown variables.

\subsection{Bernstein Polynomial}
Bernstein polynomials are parametric basis functions of degree $n$, and it is defined as\cite{Farin2002}
\begin{equation}
    b_{n,k}(s) = \binom{n}{k}s^k (1-s)^{n-k} = \frac{n!}{(n-k)!k!}\,s^k (1-s)^{n-k},
\end{equation}
where $k=0, \cdots, n$, and $s$ is the continuous parameter for $s \in [0,\,1]$.

Also, Bernstein polynomials represent the basis function of non-rational B\'ezier curves and surfaces.

\subsection{Non-rational B\'ezier Curves}
A non-rational B\'ezier curve of degree $n$ is defined as\cite{Farin2002}
\begin{equation}\label{eq:nonrational}
    {\bf B}(s) = \sum_{k=0}^{n}b_{n,k}(s){\bf P}_k,
\end{equation}
where ${\bf P}_k$ are the control points.

\section{Approximated Solutions for TPBVPs}

Consider a general form of TPBVPs given by\cite{Burden2014}
\begin{equation}   \label{eq:TPBVPdef}
    \ddot{\bf x}(t) = {\bm f}(t,\,{\bf x}(t),\, \dot{\bf x}(t)), \, \text{for }t \in [t_\text{i} ,\, t_{\text{f}}]\, \text{ with }\,{\bf x}(t_\text{i}) = {\bf x}_\text{i} \text{ and } {\bf x}(t_{\text{f}})={\bf x}_{\text{f}},
\end{equation}
where ${\bf x}(t)\in \mathbb{R}^m$ is the state vector, ${\bm f} \in \mathbb{R}^m$ is the nonlinear function, $t$ is the time, and the subscripts $\text{i}$ and $\text{f}$ indicate the initial and the final parameters, respectively. The given equation is rewritten as
\begin{equation}\label{eq:residual}
    {\bm g}(t,{\bf x}(t), \dot{\bf x}(t),\ddot{\bf x}(t))=\ddot{\bf x}(t) -{\bm f}(t,{\bf x}(t), \dot{\bf x}(t)),
\end{equation}
where ${\bm g}(t,{\bf x}(t), \dot{\bf x}(t),\ddot{\bf x}(t))\in \mathbb{R}^m$ is the residual based on the given TPBVP.

Using the definition of non-rational B\'ezier curves with order $n$ in Eq. \eqref{eq:nonrational}, $t$ and ${\bf x}(t)$ are transformed into $t(s)$ and ${\bf x}(s)$. The first and second time derivatives of ${\bf x}(s)$ are expressed as\cite{Kim2014}
\begin{equation}\label{eq:derivativex}
\begin{split}
    \dot{\bf x}(s) & = \frac{\mathrm{d}{\bf x}}{\mathrm{d}t} = \frac{\mathrm{d}{\bf x}}{\mathrm{d}s}\frac{\mathrm{d}s}{\mathrm{d}t}= \frac{\mathrm{d}{\bf x}}{\mathrm{d}s}\left( \frac{\mathrm{d}t}{\mathrm{d}s} \right)^{-1} = \frac{{\bf x}'}{t'},\\
    \ddot{\bf x}(s) & = \frac{\mathrm{d}}{\mathrm{d}t}\left( \frac{{\bf x}'}{t'} \right) = \frac{\mathrm{d}}{\mathrm{d}s}\left( \frac{{\bf x}'}{t'} \right) \frac{\mathrm{d}s}{\mathrm{d}t}= \frac{\mathrm{d}}{\mathrm{d}s}\left( \frac{{\bf x}'}{t'} \right)\frac{1}{t'} = \frac{{\bf x}'' t' - t'' {\bf x}'}{t'^3},
\end{split}
\end{equation}
where the superscripts $'$ and $''$ mean the first and second derivatives with respect to the variable $s$, respectively. 
% {\color{red}The control points ${\bf P}_k$ in Eq. \eqref{eq:nonrational} correspond to $t_k$ and ${\bf x}_k$ in Eq. \eqref{eq:derivativex}. Also, although each $t_k$ has one dimension, each ${\bf x}_k$ has $m$ dimensions.} 
The endpoints, $t_{0}$, $t_{n}$, ${\bf x}_{0}$, and ${\bf x}_{n}$, are replaced with the boundary conditions, $t_{\text{i}}$, $t_{\text{f}}$, ${\bf x}_{\text{i}}$, and ${\bf x}_{\text{f}}$, respectively because the boundary conditions of the given problem are mapped into $s\in [0,\,1]$. Using the expression provided by Eq. \eqref{eq:derivativex}, the residual function in Eq. \eqref{eq:residual} is transformed into the B\'ezier function form as
\begin{equation}\label{eq:residualBezierform}
\begin{split}
    {\bm g}_s(s,t_1, \cdots , t_{n-1},{\bf x}_1, \cdots, {\bf x}_{n-1}) & = \frac{{\bf x}'' t' - t'' {\bf x}'}{t'^3} -{\bm f}_s(t(s),{\bf x}(s), \dot{\bf x}(s)),
\end{split}
\end{equation}
where $t_1, \cdots , t_{n-1}$ and ${\bf x}_1, \cdots, {\bf x}_{n-1}$ are the intermediate control points as unknown and ${\bm f}_s$ is the transformed function in the form of B\'ezier curves. To obtain the unknown control points, the function $L$ is defined as 
\begin{equation}\label{eq:L}
\begin{split}
    L(t_1, \cdots , t_{n-1},{\bf x}_1, \cdots, {\bf x}_{n-1}) & = \int_0^1 {\bm g}_s^T {\bm g}_s\, \mathrm{d}s.
\end{split}
\end{equation}
Once the control points to minimize the function $L$ are obtained, the approximated solution in terms of B\'ezier curves for the given problem is calculated by substituting the selected control points to B\'ezier curves, and the initial value at $s=0$, $\dot{\bf x}(s)$, is determined.

\section{Numerical Studies for TPBVPs}

Before applying the method proposed to the orbit determination problem, a nonlinear TPBVP with one dimension is considered.

\subsection{Nonlinear TPBVP with One Dimension}
Consider a nonlinear TPBVP given by\cite{Burden2014}
\begin{equation}
    \ddot{x}(t) = \frac{1}{8} (32 + 2t^3 - x(t) \dot{x}(t)),\, \text{for } t\in [1,\,3] \text{ with }x(1)=17\text{ and } x(3)=\frac{43}{3},
\end{equation}
The exact solution of the given problem is $x(t) = t^2 + {16}/{t}$, and its time derivative is $\dot{x}(t) = 2t -{16}/{t^2}$.
% The residual \textcolor{blue}{function} for the given TPBVP is expressed as
% \begin{equation}
%     {g}(t,{x}(t),\ddot{x}(t)) = \ddot{x}(t) - \frac{1}{8}(32 + 2t^2 - x(t) \dot{x}(t)).
% \end{equation}
To transform the given problem into the B\'ezier function form, Eqs. \eqref{eq:derivativex} and \eqref{eq:residualBezierform} are utilized. In this work, the quadratic non-rational B\'ezier curves are considered to describe nonlinear functions. B\'ezier functions of $t$ and $x(s)$ are defined as
\begin{equation}\label{eq:tsxs}
\begin{split}
    t(s) & = (1-s)^2 t_0 + 2s(1-s) t_1 + s^2 t_2, \\
    x(s) & = (1-s)^2 x_0 + 2s(1-s) x_1 + s^2 x_2,
\end{split}
\end{equation}
% where $t_k$ and $x_k$ are the control points and $k\,(k = 0, 1, 2)$ is the index. Since the order of B\'ezier curves is 2, there are 6 control points. 
and the derivatives of $t(s)$ and $x(s)$ with respect to $s$ are found as
\begin{equation}\label{eq:deriv}
\begin{split}
    t'(s) & = 2(s-1)t_0 + 2(1-2s)t_1 + 2st_2,\\
    x'(s) & = 2(s-1)x_0 + 2(1-2s)x_1 + 2sx_2,\\
    t''(s) & = 2(t_0 - 2t_1 + t_2), \\
    x''(s) & = 2(x_0 -2x_1 + x_2).
\end{split}
\end{equation}

The residual function with B\'ezier curves is expressed as
\begin{equation}\label{eq:1dex}
    {g}_s(s,t_1,x_1) = \frac{x''t'-t''x'}{t'^3}-\frac{1}{8}\left(32+2t^3(s)-x(s)\frac{x'}{t'}\right),
\end{equation}
where $t(s)$, $x(s)$, and its derivatives are in Eqs. \eqref{eq:tsxs} and \eqref{eq:deriv}. The endpoints are $t_0=1,\,t_2=3, x_0=17$, and $x_2={43}/{3}$, and the unknown intermediate control points are $t_1$ and $x_1$. Note that the result of Eq. \eqref{eq:1dex} is not described here since the final equation with B\'ezier curves is massive.

To find the proper control points, it needs to solve the following problem:
\begin{equation}
    \text{min. }L(t_1, x_1) = \int_0^1 { g}_s^2\, \mathrm{d} s.
\end{equation}

After the control points to minimize the $L$ function are obtained, the initial value, $\dot{x}(s)$, as an initial guess is calculated by using Eq. \eqref{eq:derivativex} at $s=0$. Using the initial guess, the solution of the TPBVP is obtained using the shooting method. 

% For the numerical simulation, MATLAB is used in this work. 
The control points are obtained by using MATLAB internal function, and the obtained control points to minimize the residual function are $t_1 = 1.5833$ and $x_1=8.4757$. Using these values, the initial guess $\dot{x}(s)$ is calculated by using Eqs. \eqref{eq:derivativex} and \eqref{eq:deriv} at $s=0$, and the found initial guess is -14.6147. Since the initial condition of the exact solution is -14.0, the error between the initial guess found and exact solution is about 4.39$\%$.

\begin{figure}[h]
    \centering
    \includegraphics[width=0.58\textwidth]{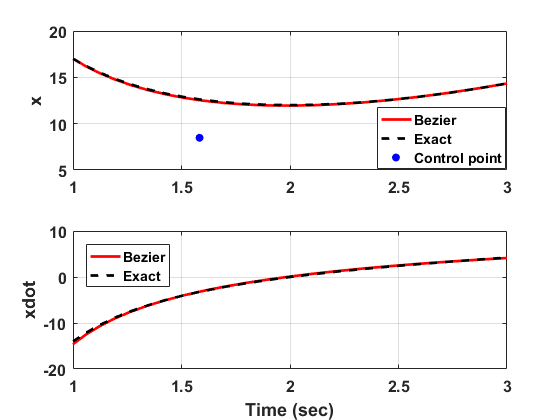}
    % \caption{Approximated and exact solution}
    % \label{fig:NLsolutionB}
% \end{subfigure}
\caption{Comparison results for the B\'ezier curves based solution and the exact solution}
\label{fig:ExactBezier}
\end{figure}

The time trajectory results for the exact solution and the one obtained by using B\'ezier curves are displayed in Fig. \ref{fig:ExactBezier}, and it is shown that the approximated solution obtained by B\'ezier curves are very similar to the exact solution.
% Also, the trajectories of the approximate solution and error based on B\'ezier curves are shown in Fig. \ref{fig:NLsolutionB} and \ref{fig:NLerrorB}.

To highlight the performance of the proposed method, the number of iterations and computational time are used. In the case of the general shooting method, an arbitrary initial guess is selected whereas the proposed method uses the initial guess obtained by using B\'ezier curves. The comparison results are listed in Table \ref{tab:compTPBVP}. 
\begin{table}[h]
  \centering
  \caption{Comparison results between the general and proposed method}
    \label{tab:compTPBVP}
    \begin{tabular}{|c|c|c|c|c|}
    \hline
    \multicolumn{3}{|c|}{Item} & General & Proposed \\ \hline
    \multicolumn{1}{|c|}{\multirow{3}[6]{*}{\shortstack{Number of\\iterations}}} & \multirow{2}[3]{*}{Shooting} & Interval & 19    & 3 \\ \cline{3-5}
     &       & Root  & 6     & 6 \\ \cline{2-5}
     & \multicolumn{2}{c|}{Total} & 25    & 9 \\ \hline 
    \multicolumn{1}{|c|}{\multirow{3}[6]{*}{\shortstack{Computational\\time (sec)}}} & \multicolumn{2}{c|}{B\'ezier} &     -  & 0.0321 \\ \cline{2-5}
     & \multicolumn{2}{c|}{Shooting} & 0.8590 & 0.2029 \\ \cline{2-5}
     & \multicolumn{2}{c|}{Total} & 0.8590 & 0.2350 \\ \hline
    \multicolumn{3}{|c|}{Improvement (\%)} & \multicolumn{2}{|c|}{72.64}  \\ \hline
    \end{tabular}%
\end{table}%

For the shooting method, the MATLAB internal function, {\it{fzero}}, is utilized in this work. This algorithm is based on Bisection method to find the root, and Bisection method needs two initial value as an interval, which is the range that there exist the root. For this reason, two kinds of iterations are performed: the one is for finding the interval that there exist the root, and the other is for finding the root using the obtained interval. Because of this process, the number of iterations is composed of ``Interval" and ``Root" in Table \ref{tab:compTPBVP}. Since the proposed method provides a similar value with the true initial condition, the number of iterations to find the interval is smaller than the one need for the general shooting method. In the case of the computational time, the proposed method has one more process to find the initial guess using B\'ezier curves. With this initial guess, the solution is found by applying the shooting method. Even though the proposed method requires one more step, the total computational time is reduced by about 72.64\%.

\subsection{Nonlinear TPBVP with Multi-Dimension: Orbit Determination}
In this subsection, B\'ezier curves are applied to a nonlinear TPBVP for the orbit determination, which is called the two-body problem. Consider the two body equation given by\cite{schaub2009}
\begin{equation}
    \ddot{\bf r}(t) = -\frac{\mu}{||{\bf r}(t)||^3}{\bf r}(t), \, \text{for }t \in [t_\text{i} ,\, t_{\text{f}}]\, \text{ with }\,{\bf r}(t_\text{i}) = {\bf r}_\text{i} \text{ and } {\bf r}(t_{\text{f}})={\bf r}_{\text{f}},
\end{equation}
where ${\bf r}(t)\in \mathbb{R}^3$ is the position vector, $\mu  \in \mathbb{R}$ is the gravitational coefficient, which is $\mu=398600$ km$^3$/s$^2$, $||{\bf r}(t)||$ is the magnitude of the position vector, $t_{\text{i}}$ and $t_{\text{f}}$ are the initial and final time, and ${\bf r}_{\text{i}}$ and ${\bf r}_{\text{f}}$ are the initial and final position vector, respectively. 

The residual function is defined as
\begin{equation}
    {\bm g}(t,{\bf r}(t),\dot{\bf r}(t),\ddot{\bf r}(t)) = \ddot{\bf r}(t) + \frac{\mu}{||{\bf r}(t)||^3}{\bf r}(t).
\end{equation}

To transform the given equation into the B\'ezier function form, the quadratic non-rational B\'ezier curves are also considered. For $t(s)$, the first equation in Eq. \eqref{eq:tsxs} is used, and $t'(s)$ and $t''(s)$ in Eq. \eqref{eq:deriv} are also considered. ${\bf r}(s)\in\mathbb{R}^3$ is defined as
\begin{equation}
\begin{split}
    {\bf r}(s) & = (1-s)^2 {\bf r}_0 + 2s(1-s) {\bf r}_1 + s^2 {\bf r}_2,
\end{split}
\end{equation}
and its derivatives with respect to $s$ are 
\begin{equation}\label{eq:derivr}
\begin{split}
    {\bf r}'(s) & = 2(s-1){\bf r}_0 + 2(1-2s){\bf r}_1 + 2s{\bf r}_2,\\
    {\bf r}''(s) & = ({\bf r}_0 -2{\bf r}_1 + {\bf r}_2).
\end{split}
\end{equation}

The residual function is transformed as
\begin{equation}\label{eq:residualTwobodyBezier}
    {\bm g}_s(s,t_1,{\bf r}_1) =\frac{{\bf r}'' t' - t'' {\bf r}'}{t'^3}  + \frac{\mu}{||{\bf r}(s)||^3} {\bf r}(s),
\end{equation}
In this problem, the unknown intermediate control points are $t_1 \in \mathbb{R}$ and ${\bf r}_1 \in \mathbb{R}^3$. Equation \eqref{eq:residualTwobodyBezier} contains the endpoints, $t_0,\,t_2, {\bf r}_0$, and ${\bf r}_2$, and these are replaced with the boundary conditions of the given problem.

To find the control points, the function $L(t_1,{\bf r}_1)$ is defined as:
\begin{equation}\label{eq:twobodyL}
    L(t_1, {\bf r}_1) = \int_0^1 {\bm g}_s ^T {\bm g}_s\, \mathrm{d}s.
\end{equation}

When the control points are selected by minimizing the $L$ function, the initial values, $\dot{\bf r}(s)$, are calculated using Eq. \eqref{eq:derivativex} at $s=0$. Using the initial guess calculated by B\'ezier curves, the solution of the two-body problem is obtained.

In this simulation, three kinds of orbits which have different eccentricities are selected. The orbits used in the simulations are circular, tundra, and Molniya orbits, and the initial conditions are assumed as near apogee, intermediate, and near perigee. In the case of the circular orbit, the orbit of the Hubble space telescope is used, and only near perigee is considered because the velocity at each position is similar. The simulation cases are listed in Table \ref{tab:cases}.

\begin{table}[htbp]
    \centering
    \caption{Simulation cases}
    \label{tab:cases}
    \begin{tabular}{|c|c|c|c|}
    \hline
    Orbit & Eccentricity & Initial position & Case \\ \hline
    Circular orbit & 0.000283 & Near perigee & 1 \\ \hline
    \multirow{3}[6]{*}{Tundra orbit} & \multirow{3}[6]{*}{0.268} & Near apogee & 2-1 \\ \cline{3-4}
    &       & Intermediate & 2-2 \\ \cline{3-4}
    &       & Near perigee & 2-3 \\    \hline
    \multirow{3}[6]{*}{Molniya orbit} & \multirow{3}[6]{*}{0.74} & Near apogee & 3-1 \\ \cline{3-4}
    &       & Intermediate & 3-2 \\ \cline{3-4}
    &       & Near perigee & 3-3 \\     \hline
    \end{tabular}
\end{table}

% Although the velocity vector in the elliptical orbit has difference direction, this unit vector is reasonable as an initial guess of the general shooting method. 

% \textcolor{red}{After the control points are selected, the initial velocity as an initial guess are also calculated by substituting the control points into the Eq. \eqref{eq:derivativex}.} 

% \textcolor{red}{In the figures, there are some errors in most of figures. In the case of case 3-3, the approximated solution using B\'ezier curves cannot follow the numerical solution especially. Since there is the rapid change of the velocity near perigee, there are large errors in the trajectory.}

\begin{table}[htbp]
  \centering
  \caption{Simulation parameters and results}
  \begin{tabular}{|c|c|c|c|c|c|c|c|c|}
  \hline
      \multicolumn{2}{|c}{\multirow{2}[4]{*}{}} & \multicolumn{3}{|c|}{Simulation parameters} & \multicolumn{4}{c|}{Results} \\
\cline{3-9}
      \multicolumn{2}{|c|}{}    & \shortstack{TOF\\(sec)} & \shortstack{${\bf r}_{\text{i}}$\\(km)} & \shortstack{${\bf r}_{\text{f}}$\\(km)} & \shortstack{$\dot{\bf r}_{\text{i}_{\text{ guess}}}$\\(km/s)} & \shortstack{$\dot{\bf r}_{\text{i}_{\text{ shooting}}}$\\(km/s)} & \shortstack{Error\\(\%)} & \shortstack{$|| \dot{\bf r}_{\text{i}_{\text{ shooting}}}||$ \\ (km/s)} \\
    \hline
    \multirow{3}[6]{*}{Case 1} & x     & \multirow{3}[6]{*}{1500} & -5641.484 & 3329.045 & 4.495 & 3.188 & 41.00 & \multirow{3}[6]{*}{7.5925} \\
\cline{2-2}\cline{4-8}          & y     &       & -3331.740 & -5754.978 & -7.470 & -6.631 & 12.65 &  \\
\cline{2-2}\cline{4-8}          & z     &       & 2204.246 & -1871.615 & -2.505 & -1.875 & 33.59 &  \\
    \hline
    \multirow{3}[6]{*}{Case 2-1} & x     & \multirow{3}[6]{*}{25000} & 15040.510 & -36285.493 & -2.614 & -2.202 & 18.70 & \multirow{3}[6]{*}{2.3829} \\
\cline{2-2}\cline{4-8}          & y     &       & 22615.098 & 13559.482 & 0.371 & 0.407 & 9.01  &  \\
\cline{2-2}\cline{4-8}          & z     &       & 45161.321 & 27077.646 & 0.740 & 0.814 & 9.02  &  \\
    \hline
    \multirow{3}[6]{*}{Case 2-2} & x     & \multirow{3}[6]{*}{15000} & -40292.402 & -17983.494 & -0.258 & -0.367 & 29.63 & \multirow{3}[6]{*}{2.9709} \\
\cline{2-2}\cline{4-8}          & y     &       & 7484.694 & -11870.227 & -1.622 & -1.320 & 22.89 &  \\
\cline{2-2}\cline{4-8}          & z     &       & 14946.572 & -23704.307 & -3.239 & -2.636 & 22.88 &  \\
    \hline
    \multirow{3}[6]{*}{Case 2-3} & x     & \multirow{3}[6]{*}{17000} & -24501.896 & 33647.418 & 3.926 & 3.005 & 30.62 & \multirow{3}[6]{*}{3.8206} \\
\cline{2-2}\cline{4-8}          & y     &       & -9999.969 & -5531.998 & -1.181 & -1.056 & 11.82 &  \\
\cline{2-2}\cline{4-8}          & z     &       & -19969.490 & -11047.131 & -2.359 & -2.109 & 11.83 &  \\
    \hline
    \multirow{3}[6]{*}{Case 3-1} & x     & \multirow{3}[6]{*}{18000} & 7062.077 & -16831.220 & -1.751 & -1.424 & 22.98 & \multirow{3}[6]{*}{1.6894} \\
\cline{2-2}\cline{4-8}          & y     &       & 19756.303 & 12838.490 & 0.358 & 0.407 & 12.16 &  \\
\cline{2-2}\cline{4-8}          & z     &       & 39452.426 & 25637.872 & 0.714 & 0.813 & 12.17 &  \\
    \hline
    \multirow{3}[6]{*}{Case 3-2} & x     & \multirow{3}[6]{*}{5000} & -17436.334 & -14505.515 & -0.715 & -0.498 & 43.64 & \multirow{3}[6]{*}{3.2791} \\
\cline{2-2}\cline{4-8}          & y     &       & 11461.543 & 1846.234 & -1.598 & -1.451 & 10.12 &  \\
\cline{2-2}\cline{4-8}          & z     &       & 22888.172 & 3686.843 & -3.191 & -2.898 & 10.12 &  \\
    \hline
    \multirow{3}[6]{*}{Case 3-3} & x     & \multirow{3}[6]{*}{5000} & -3653.531 & 17638.454 & 7.836 & 9.250 & 15.28 & \multirow{3}[6]{*}{9.6848} \\
\cline{2-2}\cline{4-8}          & y     &       & -2844.545 & 6821.862 & -2.445 & -1.285 & 90.21 &  \\
\cline{2-2}\cline{4-8}          & z     &       & -5680.425 & 13622.943 & -4.882 & -2.567 & 90.21 &  \\
    \hline
    \end{tabular}%
  \label{tab:initialguess}%
\end{table}%

Table \ref{tab:initialguess} shows the simulation conditions for each case, the calculated initial guesses using B\'ezier curves, and its errors with respect to the initial condition to be found. To find the control points for B\'ezier curves, the mid-point between the initial and final position is selected as an initial guess. The TOF, ${\bf r}_{\text{i}}$, and ${\bf r}_{\text{f}}$ are used as the simulation parameters, and $\dot{\bf r}_{\text{i}_{\text{ guess}}}$ and $\dot{\bf r}_{\text{i}_{\text{ shooting}}}$ are the initial guess obtained from B\'ezier curves and the initial value calculated by the shooting method using the obtained initial guess, respectively. Case 1 shows large errors, and the magnitude of the norm of the initial condition is also large value. In case 2 and 3, the large errors are especially observed in case 2-3 and 3-3, near perigee cases. For each case, the results seem to have a tendency that the large velocity gives the large error.

\begin{figure}[htbp]
    \centering
    \includegraphics[width = 0.58\textwidth]{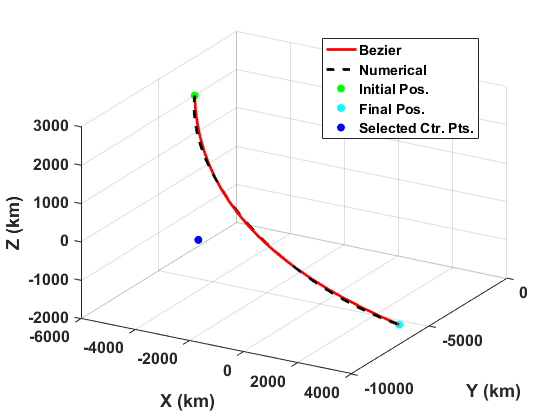}
    \caption{Approximated solution using B\'ezier curves for case 1}
    \label{fig:hubble}
\end{figure}

\begin{figure}[htbp]
\begin{subfigure}{0.32\columnwidth}
    \centering
    \includegraphics[width=1\linewidth]{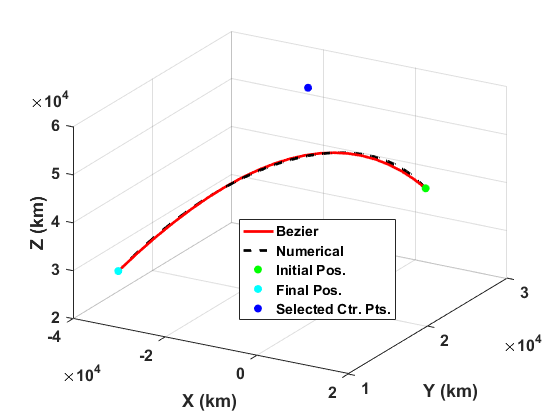}
    \caption{case 2-1}
\end{subfigure}
\begin{subfigure}{0.32\columnwidth}
    \centering
    \includegraphics[width=1\linewidth]{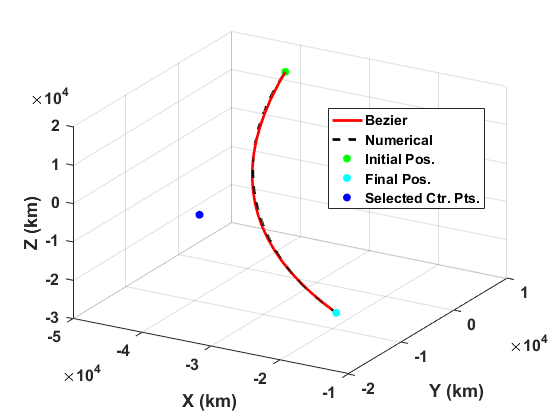}
    \caption{case 2-2}
\end{subfigure}
\begin{subfigure}{0.32\columnwidth}
    \centering
    \includegraphics[width=1\linewidth]{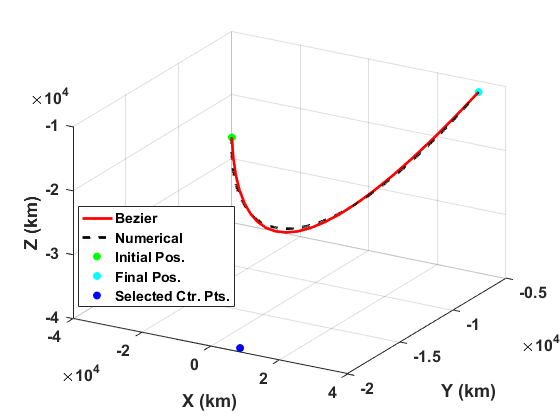}
    \caption{case 2-3}
\end{subfigure}
\caption{Approximated solution using B\'ezier curves for case 2}
\label{fig:tundra}
\end{figure}

\begin{figure}[htbp]
\begin{subfigure}{0.32\columnwidth}
    \centering
    \includegraphics[width=1\linewidth]{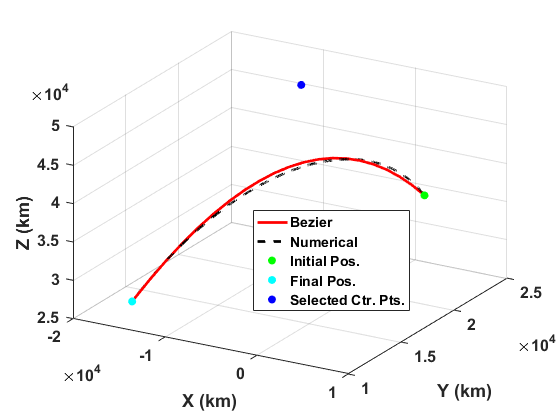}
    \caption{case 3-1}
\end{subfigure}
\begin{subfigure}{0.32\columnwidth}
    \centering
    \includegraphics[width=1\linewidth]{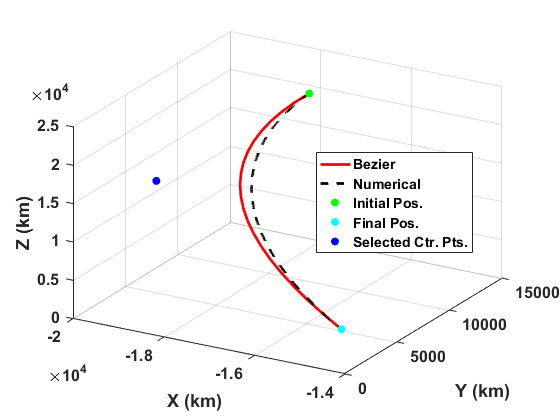}
    \caption{case 3-2}
\end{subfigure}
\begin{subfigure}{0.32\columnwidth}
    \centering
    \includegraphics[width=1\linewidth]{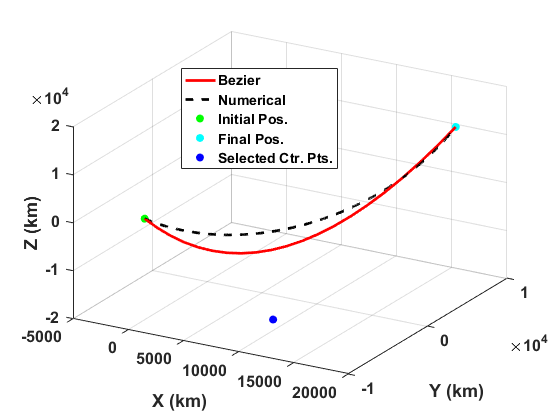}
    \caption{case 3-3}
\end{subfigure}
\caption{Approximated solution using B\'ezier curves for case 3}
\label{fig:molniya}
\end{figure}

To investigate the trajectories between the initial and final position, the approximated solution obtained by B\'ezier curves with the found control points and the numerical integration result using the initial value determined by the shooting method are displayed in Figs. \ref{fig:hubble} - \ref{fig:molniya}. The trajectories of the approximated solution look similar to the numerical integration results in most cases although the case 3-3 has the trajectory error. Since the approximated solution is used to find the initial guess to solve the given TPBVP in this work, these results show that the proposed method is applicable to guess the initial value.

For the general shooting method, a vector of the cross product between the orbit normal unit vector and the initial position unit vector is used as an initial guess because this gives the rough direction of the velocity. Table \ref{tab:comparisonOrbit} shows the computational time and the number of iterations to find the solution for each case. In the general shooting method, the large number of iterations is required for case 1 and case 2-3, and it even fails to find the solution for case 3-3. It seems that these results are affected by the large difference between the initial guess for the general shooting method and the initial condition to be found. However, the proposed method provides less number of iterations than the general shooting method. Since the initial guesses obtained by B\'ezier curves have similar values with the initial conditions to be found, the proposed method provides better results. For case 3-3, only the proposed method gives the solution while the general shooting method fails to find the solution. In the aspect of the computational time, the proposed method requires less computational time than the general shooting method for all cases. In particular, for case 1, case 2-3, and case 3-3, the proposed method shows the outstanding performance when the magnitude of the initial condition is large. These results indicate that the proposed method can be utilized to the orbit determination including highly eccentric orbits as well.

\begin{table}[htbp]
  \centering
  \caption{Computational time and number of iterations for each case}
    \begin{tabular}{|c|c|c|c|c|c|c|c|c|}
    \hline
    \multirow{3}[8]{*}{Case} & \multicolumn{3}{c|}{General} & \multicolumn{4}{c|}{Proposed} & \multirow{3}[10]{*}{\shortstack{Improvement\\(\%)}} \\
\cline{2-8}          & \multicolumn{2}{c|}{Shooting} & \multicolumn{1}{c|}{\multirow{2}[7]{*}{\shortstack{Total\\time\\(sec)}}} & \multicolumn{1}{c|}{B\'ezier} & \multicolumn{2}{c|}{Shooting} & \multicolumn{1}{c|}{\multirow{2}[7]{*}{\shortstack{Total\\time\\(sec)}}} &  \\
\cline{2-3}\cline{5-7}          & \multicolumn{1}{c|}{\shortstack{Time\\(sec)}} & \multicolumn{1}{c|}{\shortstack{\# of\\iterations}} &       & \multicolumn{1}{c|}{\shortstack{Time\\(sec)}} & \multicolumn{1}{c|}{\shortstack{Time\\(sec)}} & \multicolumn{1}{c|}{\shortstack{\# of\\iterations}} &       &  \\
    \hline
    1     & 1.5835 & 18    & 1.5835 & 0.5493 & 0.3558 & 5     & 0.9051 & 42.84   \\
    \hline
    2-1   & 5.8253 & 7     & 5.8253 & 0.4084 & 3.4314 & 4     & 3.8398 & 34.08  \\
    \hline
    2-2   & 2.8471 & 5     & 2.8471 & 0.4608 & 2.0769 & 4     & 2.5377 & 10.87   \\
    \hline
    2-3   & 5.2654 & 10    & 5.2654 & 0.5380 & 2.4071 & 4     & 2.9451 & 44.07  \\
    \hline
    3-1   & 3.7242 & 6     & 3.7242 & 0.4694 & 2.4571 & 4     & 2.9265 & 21.42   \\
    \hline
    3-2   & 1.2193 & 5     & 1.2193 & 0.3496 & 0.7598 & 4     & 1.1094 & 9.01   \\
    \hline
    3-3   & -     & -     & -     & 0.5729 & 1.3177 & 7     & 1.8906 & - \\
    \hline
    \end{tabular}%
  \label{tab:comparisonOrbit}%
\end{table}%

\section{Conclusion}
In this work, an effective method to find the initial guess of TPBVPs is proposed. The given TPBVP is transformed into B\'ezier curves form. After the control points are obtained by solving the system of nonlinear equations, the approximated solutions are calculated by substituting the control points into B\'ezier curves. Using the approximated solution, the initial guess for solving TPBVP is simply calculated, and the solution of TPBVP using the calculated initial guess is obtained through the shooting method. To validate the performance of the proposed method, numerical simulations are performed with several orbit cases. As a result, initial guesses are well found, and these initial guesses reduce the computational burden to find the solution. Also, it is shown that the solution for the highly elliptical orbit is obtained by using the method proposed only.

\bibliographystyle{AAS_publication}   % Number the references.
\bibliography{references}   % Use references.bib to resolve the labels.

\end{document}